\documentclass[12pt]{article}

\font\tenmsb=msbm10
\font\sevenmsb=msbm7
\font\fivemsb=msbm5

\newfam\msbfam
\textfont\msbfam=\tenmsb
\scriptfont\msbfam=\sevenmsb
\scriptscriptfont\msbfam=\fivemsb
\def\Bbb#1{{\fam\msbfam #1}}

\makeatletter
\ifnum\@ptsize=0  \fi \ifnum\@ptsize=2  \fi 

\newcommand\qed{{\hspace*{\fill}Q.E.D.\vskip12pt plus 1pt}}

\newcommand\sE{{\cal E}}
\newcommand\sF{{\cal F}}
\newcommand\sG{{\cal G}}

\newcommand\sI{{\cal I}}

\newcommand\sL{{\cal L}}
\newcommand\sO{{\cal O}}

\newcommand\sV{{\cal V}}

\newcommand\zed{{\Bbb Z}}
\newcommand\comp{{\Bbb C}}

\newcommand\bP{{\Bbb P}}

\newcommand\proof{{\noindent\bf Proof.\ }}

\newtheorem{theorem}{Theorem}[section]
\newtheorem{lemma}[theorem]{Lemma}
\newtheorem{corollary}[theorem]{Corollary}
\newtheorem{proposition}[theorem]{Proposition}

\newtheorem{re}[theorem]{Remark}

\newtheorem{definition}[theorem]{Definition}

\newtheorem{problem}[theorem]{Problem}
\newtheorem{note}[theorem]{Note}
\newtheorem{example}[theorem]{Example}
\newtheorem{notation}[theorem]{Notation}

\newenvironment{remark}{\begin{re}\em}{\end{re}}

\begin{document}
\title {Ample Vector Bundles and Branched Coverings, II} \author{Thomas Peternell\and Andrew J. Sommese}
\date{november 5, 2002}

\maketitle

\vspace*{-0.5in}\section*{Introduction}

To a degree $d$ covering $f: X \to Y$ of projective manifolds
there is naturally attached a holomorphic rank $(d-1)-$vector
bundle $\sE_f,$ just by dualizing $f_*(\sO_X)$ and dividing by the
trivial factor. This bundle reflects interesting properties of $f$
or $X.$ It has certain positivity properties: it is nef on the
generic curve, as observed in [PS00], the first part to this
paper, and even more, it is nef on every curve not contained in
the branch locus of $f$ (Lazarsfeld). Often it has a lot of
sections ([PS00]). For special $Y$ much more is true: if $Y =
\bP_n$, then $\sE$ is ample and spanned (Lazarsfeld); the same
being true for some rational homogeneous manifolds (Kim, Manivel)
and conjecturally for all (as long as $b_2 = 1$).
 In the first part to this paper we investigated Fano manifolds
and proved nefness in certain cases (del Pezzo manifolds of higher
degrees). We suspected
that we might have ampleness as well. \\
 The first result in this paper is that $\sE$ is not always ample
(2.1): we provide a triple covering over the Fano threefold $V_5$,
appearing as a linear section of the Grassmannian $G(2,5) \subset
\bP_9,$ for which $\sE$ is not ample. However the only curves on
which $\sE$ fails to be ample, are lines (not necessarily in the
branch locus). Then (2.2) we give a criterion for ampleness, and
as a consequence we see that given a
covering $X \to V_5,$ then $\sE$ is ample if and
only if $\sE$ is ample on every line. Moreover we show that $\sE$
is always $1-$ample. \\
The second part of the paper is devoted to the study of low degree
coverings (of Fano manifolds, usually). Again, a result of
Lazarsfeld is the starting point: if $f: X \to \bP_n $ has degree
$d$, then the restriction maps
$$ H^i(\bP_n) \to H^i(X) $$
are surjective for $i \leq n+1-d.$  We give various
generalizations, e.g., we prove the following: Let $f: X \to Y$
have degree 3. If $Y$ is simply connected and $h^{1,1}(Y) = 1,$
then $q(X) = q(Y) = 0.$ Actually, even if $Y$ is not simply
connected, we have
$q(X) = q(Y)$ unless $f$ is \'etale.\\
 We also show that if say $\sE$ is spanned and the map defined by
the global sections of $\sO(1)$ on $\bP(\sE) $ is ``semi-small'',
then we have the same topological conclusion as if $\sE$ were
ample. The decisive role here is played by the top Chern class
$c_{d-1}(\sE),$ in connection with the hard Lefschetz theorem. It
is an interesting question whether the positivity of
$c_{d-1}(\sE)$, in a differential-geometric or algebraic sense,
suffices to guarantee the hard Lefschetz theorem (cupping with
$c_{d-1}(\sE)).$ \\
 In sect.4 we apply the positivity of $\sE_f$ to study the
ramification divisor $R \subset X.$ E.g. we prove that $R$ is
ample if $Y$ is a Grassmannian, the cases of projective space and
quadrics (of dimension at least 3) being known since the 80's. We
also show that the ramification divisor of a covering over the
Fano threefold $V_5$ is always ample unless the covering arises by
base change from the Example 2.1 (in that case $R$ is nef but not
ample).

\tableofcontents

\section{Preliminaries}

\begin{notation}{\rm Let $X$ and $Y$ be projective manifolds of dimension $n$ and
let $f: X \to Y$ be a finite covering of degree $d$. Then we set
$$ \sE = [f_*(\sO_X)/\sO_Y]^*,$$
a vector bundle of rank $d-1$ on $Y,$ {\it the vector bundle associated to $f.$}}

\end{notation}

For details and many properties we refer to [PS00] and the
references given there. The bundle $\sE$ has certain positivity
properties which depend on the geometry of $Y.$ For the
convenience of the reader here are some results relevant for this
paper.

\begin{theorem} Let $f: X \to Y$ be a covering of projective manifolds.
\begin{enumerate}
\item ([La00]) $\sE$ is nef modulo the branch locus $B \subset Y,$ i.e., $\sE \vert C$ is nef for
every curve $C \not \subset B.$
\item ([PS00]) Let $Y$ be a del Pezzo manifold of degree at least $5.$ Then $\sE$ is spanned by global sections,
in particular $\sE$ is nef.
\end{enumerate}
\end{theorem}

\begin{notation}{\rm Let $X$ be a projective manifold and $\sL$ a line bundle on $X.$
$\sL$ is semi-ample, if some positive multiple $\sL^{\otimes m}$ is spanned by global sections.
The Stein factorization $\phi: X \to Z$ of the map associated to the linear system of
$\sL^m$ is called the {\it map associated to} $\sL$. So by definition, $Z$ is normal and $\phi$ has connected
fibers. }
\end{notation}

\begin{notation}{\rm A line $L$ is said to be $k-$ample, if $L$ is semi-ample and the associated map has at most $k-$dimensional
fibers. \\
A vector bundle $\sE$ is $k-$ample, if $\sO_{\bP(\sE)}(1) $ is $k-$ample on $\bP(\sE).$ }
\end{notation}

\begin{notation}{\rm Let $X$ be a projective manifold. Then $H^k(X) $ will
always denote cohomology with complex coefficients:
$H^k(X) = H^k(X,\comp).$ }
\end{notation}

\section{Coverings of the del Pezzo $3-$fold $V_5$.}

In the paper \cite{PS00} we suspected that the vector bundle $\sE$
associated to a branched covering $f: X \to Y$ over a Fano
manifold $Y$ with $b_2(Y) = 1$ is always ample. Here we present an
example that this is not quite true, however the deviation from
being ample is very small.

\begin{example}{\rm (1)  Let $Y$ be a Fano 3-fold of index
$2$ and of type $V_5.$ By definition, $Y$ is a submanifold of the
Grassmannian $G(2,5) \subset \bP_9$, given by three linear sections. If
$H$ is the ample generator of ${\rm Pic}(Y) = \zed,$ then $H^3 =
5.$ We consider the family of lines in $Y,$ the geometry of which
being studied in \cite{FN89}. The parameter space $T$ of this
family is just the projective plane, so that the graph $p: G \to
Y$ is a $\bP_1-$bundle $q: G \to \bP_2 = T.$ By \cite{FN89}, $p$
is finite of degree $3$; this will be our example: set $X = G$ and
$f = p.$ Let $R \subset X$ be the branch locus of $f;$ then
 $$
 c_1(\sE) = {{f_*(R)} \over {2}}
 $$
by \cite{PS00}(1.5). Since the
ramification divisor $B \subset Y$ is an element of $\vert -K_Y
\vert,$ we conclude from \cite{FN89}(2.3) that $f_*(R) = B,$ hence
  $$ c_1(\sE) = H.$$
Since $\sE$ is spanned by \cite{PS00}, the rank $2$-bundle $\sE$
splits on every line $l \subset Y$ as $$ \sE = \sO \oplus
\sO(1).$$ Hence $\sE$ cannot be ample. \\
 (2) Let us have a further look at the geometry of the bundle $\sE.$
First notice that because of the nefness of $\sE,$ the manifold $\bP(\sE) $
is Fano. Let $$ \psi: \bP(\sE) \to
Z$$ be the second extremal contraction of $\bP(\sE)$ besides the
projection $\pi: \bP(\sE) \to Y.$ Then $\psi$ is given - up to Stein
factorization - by the linear system $\vert \sO_{\bP(\sE)}(1)
\vert.$ Being careful, let $\phi: \bP(\sE) \to Z'$ be the map
defined by $\vert \sO_{\bP(\sE)}(1) \vert.$ \\
 In order to analyze $\phi,$ we compute $c_2(\sE).$
Notice first that $H^4(Y,\zed)$ is generated by $[l].$ Let $W \in
\vert H \vert$ be general and set $\tilde W = f^*(W). $ Since $W
\cdot l = 1,$ the smooth surface $\tilde W$ is generically a
section of $q: X \to T,$ hence $\tilde W$ (and $W$) are rational.
Let $g = f \vert \tilde W.$ \\
 By \cite{PS00}(1.15) we have $$ c_2(\sE_g) = 3 \chi(\sO_{W}) -
\chi(\sO_{\tilde W}) + {{(K_W + H_W)} \over {2}} \cdot H_W^2,$$
hence $c_2(\sE_g) = 2.$ Thus $c_2(\sE) \cdot H = c_2(\sE_g) = 2,$
and consequently $$ c_2(\sE) = 2l.$$ Having in mind $c_1^2(\sE) =
5l, $ we obtain $$ c_1(\sO_{\bP(\sE)}(1))^4 = c_1^3(\sE) - 2
c_1(\sE)c_2(\sE) = 1 > 0.$$ Thus $\sO_{\bP(\sE)}(1)$ is big, and
therefore the map $\phi$ is generically finite. \\ Now
Riemann-Roch gives $\chi(Y,\sE) = 5.$ Since $H^q(Y,\sE) = 0$ for
$q \geq 1,$ \cite{PS00}, we have $h^0(\sE) = 5$. Therefore $Z'
\subset \bP_4$ and since $\phi$ is generically finite, we have $Z'
= \bP_4.$ Moreover by
$c_1(\sO_{\bP(\sE)}(1))^4 = 1$, we see that
$\phi$ is generically $1:1$. Zariski's Main Theorem now
implies that $\phi = \psi$ and $Z = Z'.$ \\
Since $Z$ is smooth, $\psi$ cannot be a small contraction; so let
$D\subset \bP(\sE)$ be the irreducible divisor contracted
by the Mori contraction $\psi.$ We are going to show that $\psi (D)$ is a
2-dimensional
submanifold of $\bP_4 = Z$ and $\psi$ is the blow-up of $\psi
(D).$ We first show that $\dim \psi (D) = 2.$ In fact, if $\dim \psi
(D) = 0,$ then $Y$ would be dominated by $\bP_3$ and if $\dim \psi
(D) = 1,$ then $Y$ would be covered by linear images of $\bP_2$,
which is absurd. Thus $\dim \psi(D) = 2.$ In order to show that $\psi(D)$ is smooth
and $\psi$ is the blow-up, it is - by Ando's theorem \cite{An85} - sufficient to
verify that $\psi \vert D$ is equidimensional, i.e.,  has only 1-dimensional fibers.
Notice that the exceptional sections in $\bP(\sE \vert
l) = \sO_l \oplus \sO_l(1)$, $l$ any line in $Y$, are contracted, so that every line
$l \subset Y$ has a unique lifting $l_0$ which is contracted by $\psi.$
Therefore the map $X \to Y$ lifts to a map $X \to \bP(\sE)$ whose image is just
$D.$ Then the projection $X \to \bP_2$ yields a map $D \to \bP_2$ which must be
$\psi \vert D.$ Therefore $\psi \vert D$ has only 1-dimensional fibers so that
Ando's theorem applies. Moreover we have seen that $\psi (D) = \bP_2$
and the exceptional sections $l_0$
are the only curves contracted by $\psi.$ Therefore the only
curves $C \subset Y$ on which $\sE$ is not ample are lines. We
also note that $\sE$ is $1-$ample. \\
 (3) As already mentioned, $\chi(Y,\sE) = 5.$ Since
$H^q(Y,\sE) = 0$ for $q \geq 1,$ \cite{PS00}, we have $h^0(\sE) =
5$. Since $\sE$ is spanned, we obtain a holomorphic map $$ \tau:
Y  \to G(2,5).$$ Since $b_2(Y) = 1,$ this map is of course finite
onto its image; we are going to prove that $\tau $ is actually an
embedding and that the image is a linear section of the Pl\"ucker
embedding $G =G(2,5) \subset \bP_9.$ If $V$ denotes the
tautological rank $2-$bundle on $G$, we have $\sE \simeq \tau^*(V).$
Denoting by $L$ the ample generator of ${\rm Pic}(G),$
then we obtain by taking determinants that $H = \tau^*(L).$ Since
$H^3 = 5$ and since $L^6 = 5,$ the map $\tau $ must be generically
finite and $\tau (Y) \subset G$ is a linear section, i.e.,  $\tau
(Y) = G \cap L_1 \cap L_2 \cap L_3$, where $L_i $ are hyperplane
sections of the Pl\"ucker embedding $G \subset \bP_9.$ By
Zariski's Main Theorem, it only remains to show that $\tau (Y) $
cannot be non-normal. This is done as follows. First observe that
$\tilde S = \tau^{-1}(S) $ is a general linear section of $Y$ and therefore
a smooth del Pezzo surface with $K_{\tilde S}^2 = 5.$ On the other hand, by
the classification of non-normal del Pezzo surfaces, Reid \cite{Re95}, $\tilde S$
must be $\bP_2$ or a ruled surface. This is a contradiction.
}
\end{example}

Recall that a vector bundle $\sE$ is very ample by definition if
the line bundle $\sO_{\bP(\sE)}(1)$ is very ample. A line in a Fano manifold with $b_2 = 1$
is a rational curve $l$ such that $H \cdot l = 1$ where $H$ denotes the ample generator.
Towards positive results we first prove:

\begin{theorem} Let $f: X \to Y$ be a covering from a projective manifold
to the Fano manifold $Y$. Suppose $b_2(Y) = 1$ and that
$Y$ has index at least $2.$ Let $H$ be the ample generator of ${\rm Pic}(Y).$
Suppose furthermore that
\begin{enumerate}
\item $\sE = \sE_f$ is spanned
\item there exists $S \in \vert H \vert$ such that $\sE_S$ is very ample
\item for all lines $l \subset Y,$ the restriction $\sE_l$ is  ample.
\end{enumerate}
Then $\sE$ is ample.
\end{theorem}

\bigskip

\proof Let $\psi: \bP = \bP(\sE) \to Z$ be the map defined by
$\vert \sO_{\bP(\sE)}(1) \vert$. Assume that $\sE$ is not ample,
i.e.,   $\sO_{\bP(\sE)}(1)$ is not ample. Then there exists a curve
$C \subset \bP$ with $\dim \psi (C) = 0. $ We claim that the restriction map
$$r: H^0(\bP,\sO_{\bP(\sE)}(1)) \to H^0(\bP(\sE_S),\sO_{\bP(\sE_S)}(1))
$$ is surjective. In fact, this is guaranteed by $$ H^1(\sE(-H)) =
0.$$ Now $$ H^1(Y,\sE(-H)) = H^{n-1}(Y,\sE^* \otimes H \otimes
K_Y) = H^{n-1}(Y,f_*(\sO_X) \otimes H \otimes K_Y)$$ since $H
\otimes K_Y$ is negative by our assumption on the index. Now the
last group is just $$ H^{n-1}(X,f^*(K_Y \otimes H)) $$ which is
zero by Kodaira vanishing. \\
Using the surjectivity of $r$ and the very
ampleness of $\sE_S$, we conclude $$0 \leq  C \cdot \bP(\sE_S)
\leq 1.$$ Since $\dim \pi(C) = 1$ ($\pi$ denoting the projection
$\bP \to Y$), we conclude first $$0 \leq \pi(C) \cdot H \leq 1$$
and thus that $\pi(C) $ must be a line in $Y.$ But by the choice of
$C$, the bundle $\sE \vert \pi(C) $ is not ample, contradicting our
assumption (3).  \qed

\begin{note} {\rm
We observe the following (and restrict ourselves to dimension $2$
for simplicity). Let $S$ be a smooth projective surface and $\sE$
a vector bundle on $S.$ Then $\sE$ is ample, if the following
holds.\\
\\
(*)  For all $x,y \in S$ there exists a curve $C$ passing through
$x$ and $y$ such that $\sE \vert C$ is very ample and such that $$
H^0(S,\sE) \to H^0(C,\sE \vert C) $$ is surjective. \\

In fact, (*) first implies that $\sO_{\bP(\sE)}(1) $ is spanned; let $\phi:
\bP(\sE) \to Z$ be the associated map. We proven even that $\phi$ is $1:1,$
implying in particular ampleness of  $\sO_{\bP(\sE)}(1) ,$ hence of $\sE.$
So let $x,y \in \bP(\sE).$ We may assume that $\pi(x) \ne \pi(y),$ otherwise
clearly $\phi(x) \ne \phi(y).$ Let $C \subset X$ be a curve containing $\pi(x)$
and $\pi(y)$ such that (*) holds for $C.$ Then $\phi \vert \bP(\sE) $ is
an embedding, hence $\phi(x) \ne \phi(y).$ }

\end{note}

In the special case that $Y = V_5,$ the last theorem can be improved:

\begin{theorem} Let $f: X \to Y$ be a covering from the smooth 3-fold
$X$ to the Fano 3-fold $Y$ of type $V_5.$ Suppose that $\sE_f =
\sE$ is ample on every line $l \subset Y.$ Then $\sE$ is ample.
\end{theorem}

\proof Since $\sE$ is spanned by \cite{PS00},(4.2), it suffices by
Theorem 2.2 to prove that $\sE_S$ is very ample for the general $S
\in \vert H \vert.$ We first follow the arguments of
\cite{PS00}(3.2). Since $K_S^2 = 5,$ we can realize $S$ as a
blow-up $p: S \to \bP_2$ of the plane in $4$ points in general
position. Let $\tilde H \subset \bP_2$ be a hyperplane and set $H
= p^*(\tilde H). $ Then we can write $$ K_S = -3H + \sum_{i=1}^4
E_i.$$ Let $H_i = H - E_i.$ Then every $H_i$ is semi-ample and
defines a conic bundle structure $\phi_i: S \to \bP_1$ with
exactly $4$ singular fibers. \\ We are going to apply Note 2.3 to
the fibers and the fiber components of the maps $\phi_i.$ All
fibers resp. fiber components $C_i$ of $\phi$ are smooth rational,
whence ampleness and very ampleness is the same.  \\
(I) If $p(x) = p(y),$ we argue as
follows. There exists some $i$ such that $x,y \in E_i;$ denote for
simplicity $l = E_i.$ Then $l$ is a line in $Y,$ hence $\sE \vert
l$ is ample and thus very ample. Thus we only need to show that $$
H^0(S,\sE) \to H^0(l,\sE_l)$$ is onto. This is certainly
guaranteed by $$ H^1(S,\sE (-l)) = 0.$$ Consider the conic bundle
$\phi_i$ which contracts $l.$ Let $l'$ be the $(-1)-$curve in $S$
such that $l+l'$ is a fiber of $\phi_i.$ Then $$ H^1(\sE_S(-F)) =
0$$ since $-K_S - F$ is big and nef (apply [PS00,3.1]). By $$ 0
\to \sE(-F) \to \sE(-l) \to \sE_{l'}(-1) \to 0$$ we conclude $$
H^1(\sE(-l)) = 0,$$ settling the case $p(x) = p(y).$ \\
(II) So we
shall assume $p(x) \ne p(y).$ Let $x_i = p(E_i).$ Choose a
possible singular conic $\tilde C = \tilde C(x,y)$ passing through
$p(x),p(y),x_1,x_2$ and $x_3$ (including the infinitesimal case
$p(x) = x_i$ or $p(y) = x_j$).  Then we obtain a possible singular
curve $$ C \in \vert 2H - E_1 - E_2 - E_3 \vert$$ passing through
$x$ and $y.$ We are going to prove that $$ H^1(\sE(-C)) = 0. \eqno
(1)$$ Let $\tilde S = f^{-1}(S),$ which is a smooth surface in
$X,$ since $S$ is chosen generic. The covering map $\tilde S  \to
S$ is again called $f.$ Since $C$ is big and nef, Kawamata-Viehweg
vanishing gives $$ H^1(\sE(-C)) = H^1(\sE^* \otimes C \otimes K_S)
= $$ $$ = H^1(f_*(\sO_{\tilde S}) \otimes C \otimes K_S) =
H^1(f^*(K_S \otimes C)) = H^1(f^*(-H \otimes E_4)).$$ Now consider
the map $g: = \phi_4  \circ f: \tilde S \to \bP_1.$ We claim that
$g_4$ has connected fibers. In fact, otherwise Stein factorization
gives a factorization $\tilde S \to B \to \bP_1$ with a finite map
$B \to \bP_1.$ Hence the ramification divisor $R_S \subset S$ of
$f$ contains some fibers of $\phi_4.$ On the other hand, the
ramification divisor $R \subset Y$ of the full map $f$ is an ample
divisor. Moreover in the decomposition $R = \sum r_i R_i$ every
$R_i$ is an ample irreducible hypersurface and therefore $R_i \cap
S$ is an irreducible ample divisor for general $S,$ contradiction.
Therefore $g$ has connected fibers. \\ Let $F$ be a general fiber
of $g.$ Then $\sI_F = g^*(\sO(-1)) = f^*(-H+E_4).$ Now the Leray
spectral sequence yields $$ h^1(\tilde S,f^*(-H+E_4)) \leq
h^1(\bP_1,\sO(-1)) + h^0(\bP_1,R^1g_*(\sO_{\tilde S})).$$ Since
$R^1g_*(\sO_{\tilde S})$ is a semi-negative vector bundle, our
claim (1) follows. \\ In order to verify 2.3(1), it remains to
show that $ \sE_C$ is ample. Note that $\tilde C$ cannot be a
double line; otherwise the 3 points $x_1,x_2,x_3$ would be on a
line, which is impossible since $-K_S$ is ample. If $\tilde C =
\tilde l_1+ \tilde l_2$ is a sum of two different line, we are
already done, since $\sE \vert l_i$ is ample by assumption, hence
$\sE_{l_i}$ is very ample, and then $\sE_C$ is very ample, too. It
remains to show that $\sE_C$ is ample in case $\tilde C$ is a
smooth conic. Suppose to the contrary that $\sE_C$ is not ample.
Then by [PS00,1.3], $X_C = f^{-1}(C)$ is disconnected. For this
argument we need to observe that $C \not \subset R.$ This follows
from the above observation that otherwise $C$ would be an ample
curve on $S$ (again we use the genericity of $S$). Now $C$ is big
and nef but not ample: let $l$ be the strict transform of the line
joining say $x_1$ and $x_2$; then $C \cdot l = 0.$ \\ Now consider
the family of ``conics''  $(\tilde C_t)_{t \in T}$ in $Y.$ Say $C
= C_{t_0}.$ Then of course $f^{-1}(C_t)$ is disconnected for
general $t \in T$ and we claim that this has to be true for all
$t$. In fact, this follows, e.g., from the semi-continuity of
$h^0(\sO_{X_{C_t}})$ or from the openness of ampleness in families
(again applying [PS00,1.3]). Now consider $t_1  \in T$ such that
$C_{t_1}$ splits into $2$ lines $l_i.$ Then say $\sE_{l_1}$ is not
ample, contradicting our assumption. Therefore (*) is verified. \qed


\begin{lemma} Let $f: X \to Y $ be a covering of projective
varieties. Let $D$ be an effective nef and big divisor on $Y.$ Then $\sE_f \vert D$ does not have
a trivial direct summand.
\end{lemma}

\proof Suppose that $$\sE_f \vert D = \sO_D \oplus V.$$
Then $h^0(\sE^*_f \vert D) \ne 0$ and therefore $h^0(f_*(\sO_X) \vert D) \geq 2.$
But then
$$ h^1(f_*(\sO_X) \otimes \sO_Y(-D)) = h^1(f^*(\sO(-D)) \ne 0.$$
This contradicts the Kawamata-Viehweg vanishing theorem, $D$ being big and nef.
\qed

\begin{corollary} Let $f: X \to Y$ be a covering of smooth $n-$dimensional projective varieties.
Suppose that ${\rm Pic}(Y) = \zed$ and that $\sE$ is spanned.
Then $\sE_f$ is $(n-2)-$ample.
\end{corollary}

\proof  Let $\psi: \bP = \bP(\sE) \to Z$ be the map
defined by  $\vert \sO_{\bP(\sE)}(1) \vert$. Notice that $\psi
\vert \pi^{-1}(y)$ is biholomorphic, hence every fiber $F$ of
$\psi$ meets $\pi^{-1}(y)$ in at most one point.
Hence $\pi \vert F$ is biholomorphic onto its image. \\
Suppose
$\sE$ is not $(n-2)-$ample. Then we can find a fiber $F$ of $\psi$
with $\dim F \geq n-1.$ Choose a general irreducible variety $S
\subset F$ of $\dim S = n-1$ and let $S' = \pi(S).$ Then $S$ is a section of
$\bP(\sE_{S'}),$ given by a rank $1$-quotient $$ \sE_{S'} \to \sL
\to 0.$$
Since $\sE$ is
spanned, clearly $\sL =\sO_{S'}$ and moreover
$$ \sE_{S'} = \sO_{S'} \oplus V $$
with some bundle $V.$ Since ${\rm Pic}(Y) = \zed$, the divisor $S'$ is ample and Lemma 2.5 provides
a contradiction. \qed

\begin{corollary} Let $f: X \to Y$ be a covering from the smooth
$3-$fold $X$ to the Fano $3-$fold $Y$ of type $V_5.$ Then $\sE =
\sE_f$ is $1-$ample.
\end{corollary}

\proof $\sE$ is spanned by [PS00,4.2].  \qed

\section{Topology of low degree coverings of Fano manifolds}

In \cite{La80} R.Lazarsfeld proved a remarkable theorem on the
topology of branched covering of projective space:

\begin{theorem} Let $f: X \to \bP_n$ be a branched cover of
degree $d.$ Then the pull-back maps $$ H^i(\bP_n) \to H^i(X)$$ are
surjective for $i \leq n+1 - d.$
\end{theorem}

Here as always in this section $H^q(X)$ denotes cohomology with
values say in $\comp$. We will fix for this section a finite cover
$f: X \to Y$ of projective manifolds of degree $d$. The associated
bundle $\sE_f$ is simply denoted $\sE$. We address the question
to what extent this theorem holds more generally, say for Fano
manifolds instead of $\bP_n.$
We assume that $\pi_1(Y) = 0;$ since otherwise we can produce
counterexamples via factorization over finite \'etale covers. Also
we shall assume $h^{1,1}(Y) = 1$ in order to avoid maps, e.g., to
curves which immediately produce counterexamples, too. So we are
reduced to Fano manifolds with $b_2 = 1$, to Calabi-Yau manifolds
or symplectic manifolds with $h^{1,1} = 1$ and to simply connected
manifolds with ample canonical bundles (and $h^{1,1} = 1$).

We first note that the proof of Lazarsfeld actually gives

\begin{proposition} The canonical maps
$$ H^i(Y) \to H^i(X) $$
are surjective ($i \leq n+1-d$), if
$$ \wedge c_{d-1}(\sE) : H^{n+1-d+l}(X) \to H^{n-1+d+l}(X) $$
is surjective for all $l \geq 0.$
Actually, fixing $l \geq 0,$
$$ H^{n+1-d-l}(Y) \to H^{n-1+d-l}(X) $$
is surjective if
$$ \wedge c_{d-1}(\sE) : H^{n+1-d+l}(X) \to H^{n-1+d+l}(X) $$
is surjective.
\end{proposition}

\begin{corollary} Suppose that $\sE$ is ample. Then
$ H^{i}(Y) \to H^{i}(X)$  is surjective for $i \leq n+1-d.$ More
generally, if $\sE$ is $k-$ample, then  $ H^{i}(Y) \to H^{i}(X)$
is surjective for $i \leq n+1-d-k.$
\end{corollary}

\proof If $\sE$ is $k-$ample, then $\bigwedge c_{d-1}(\sE)$ is
surjective in the given range by \cite{So78}.  \qed

\begin{proposition} Suppose that $h^{d-1,d-1}(Y) = 1$ and
that $c_{d-1}(\sE) \ne 0.$ Then  $ H^{i}(Y) \to H^{i}(X)$  is
surjective for $i \leq n+1-d.$
\end{proposition}

\proof By the assumption $h^{d-1,d-1}(Y) = 1$ we can write $$
c_{d-1}(\sE) = c [\omega^{d-1}]$$ for a fixed K\"ahler form
$\omega$ and some real number $r.$ Now $c \ne 0$ by our assumption
and thus the Hard Lefschetz theorem gives the result via (3.2).
\qed

To obtain $H^1(Y) \to H^1(X) $ surjective, it is actually
sufficient to assume that $h^{1,1}(Y) = 1$ (and $c_{d-1}(\sE) \ne 0$ of course). To see this, reduce yourself to $d
= n$ by taking hyperplane sections and apply 3.4.

\begin{corollary} Suppose that $Y$ is a Fano $3-$fold of type $V_5.$
Then  $ H^{i}(Y) \to H^{i}(X)$  is surjective for $i \leq 4-d.$
\end{corollary}

Essentially this says that triple covers $X \to V_5$ have $q(X) =
0.$ \bigskip

\proof If $\sE$ is ample, this follows from (3.3). In general, we
only know that $\sE$ is $1-$ample, so we want to apply (3.4) and
need to show that $c_2(\sE) \ne 0$ (we are left with $d = 3$).
Suppose $c_2(\sE) = 0.$ Since $\sE$ is spanned, $\sE$ will have a
section without zeroes. So we obtain a sequence $$ 0 \to \sO_Y \to
\sE \to \sL \to 0 $$ with $\sL$ ample. By Kodaira vanishing
$H^1(\sL^*) = 0,$ hence $\sE = \sO_Y \oplus \sL.$ This contradicts
$h^0(\sE^*) = 0.$ \qed

Using stronger tools, we can prove much more:

\begin{corollary} Let $d = 3$ and $\dim Y \geq 3$. Suppose
${\rm Pic}(Y) = \zed.$ Then $q(X) = q(Y) = 0.$
\end{corollary}

\proof By taking hyperplane sections of $Y$ we reduce ourselves to
$\dim Y = 3.$ Then by (3.4) we may assume that $c_2(\sE) = 0.$
Since $\sE$ might not be spanned, the arguments of (3.5) have
to be modified. Since $c_1(\sE) > 0,$ the Bogomolov-Yau inequality
is violated and therefore $\sE$ is $H-$unstable for a
polarization $H.$ Let $\sL \subset \sE$ be the maximal
$H-$destabilizing subsheaf. Then we obtain a sequence $$ 0 \to \sL
\to \sE \to \sI_Z \otimes \sL' \to 0 \eqno (*)$$ with $Z \subset
Y$ of codimension $2$ or empty and $\sL'$ another line bundle on
$Y.$ The destabilizing property $$ c_1(\sL) \cdot H^2 > {{1} \over
{2}} c_1(\sE) \cdot H^2 $$ together with $h^{1,1}(Y) = 1$ implies
that $\sL$ is ample. Moreover $\sE$ is generically nef (see
[PS00]), thus $\sL'$ is generically nef and therefore nef. Now by
(*) $$ 0 = c_2(\sE) = c_1(\sL) \cdot c_1(\sL') + [Z],$$ hence $Z =
\emptyset$ and $c_1(\sL') = 0.$ Since ${\rm Pic}(Y) = \zed,$ the bundle
$\sL'$ is trivial. Dualizing (*) we therefore obtain a section in
$\sE^*,$ which is absurd. \qed

Note that ${\rm Pic}(Y) = \zed $ is satisfied as soon as $\pi_1(Y) = 0$
and $h^{1,1}(Y) = 1.$\\
\\
If $d \geq 4,$ then $\sE$ has rank $d-1 > 2,$ and it is much
harder to make use of the vanishing $c_{d-1}(\sE) = 0.$
Nevertheless it seems reasonable to expect that $c_{d-1}(\sE) \ne
0$ if $Y$ is Fano; and maybe much more should be true (see the end
of this section and sect.4). \\

In 3.5/3.6 flat (generic) quotients of $\sE$ appeared and led to
contradictions. This leads to the following considerations. \\ In
[PS00] (2.1 and its proof) we have seen that if $\dim Y = 1$ and
if $\sE_f$ is not ample, then $f$ factors via an \'etale cover.
This can be restated as follows: if $\sE_f$ has a numerically flat
quotient, then $f$ factors. One may ask whether this is true also
in higher dimensions.

\begin{problem} Let $f: X \to Y$ be a finite cover of projective
manifolds. Suppose that there exists an epimorphism $$ \sE_f \to
\sV \to 0 $$ with a numerically flat vector bundle $\sV.$ Does
there exist an \'etale cover $h: Z \to Y$ and a finite cover $g: X
\to Z$ such that $f = h \circ g$, together with an epimorphism
$\sE_h \to \sV \to 0?$
\end{problem}

We prove that this problem has a positive answer even in a slightly more general
context.

\begin{proposition}
Suppose that there is an epimorphism $$ \sE_f \to \sI_Z \otimes
\sV \to 0 $$ with a numerically flat vector bundle $\sV$ and an
analytic set $Z$ of codimension at least 2. Suppose that $f$ is
not \'etale. Then
\begin{enumerate}
\item for a general complete intersection curve $C \subset Y$
there exists a factorization $X_C \to Z(C) \to C$ with $h_C: Z(C)
\to C$ \'etale and an epimorphism $\sE_{h_C} \to \sV_C;$
\item $\deg f $ is not prime;
\item  $f$ factors via an \'etale map $h: Z \to Y$ such that $\sV$ is a quotient of
$\sE_h.$
\end{enumerate}
\end{proposition}

\proof (1) follows from [PS00,2.1] (including the proof) by
remarking that $X_C$ is connected (Bertini). \\
 (2) is trivial from (1). \\
 (3) First we construct $h$ outside a set $W \subset Y$ of codimension at least 2.
Fix an ample divisor $H$ on $Y.$ Let $C \subset Y$
 be a general complete intersection curve cut out by $n-1$ elements in $\vert mH \vert,$ where
$m$ is sufficiently large. Let $\sG$ be the kernel of $\sE \to
\sI_Z \otimes \sV.$ By our assumption, $\sG_C$ is a
$H_C$-destabilizing subsheaf of $\sE_C.$ Let $\sF_C$ be the
corresponding maximal destabilizing subsheaf; then by
Mehta-Ramanathan, $\sF_C$ extends to a torsion free coherent
subsheaf $\sF \subset \sE;$ we may assume $\sF \subset  \sE$ is
saturated. Let $T = \sE/\sF,$ a torsion free sheaf. By throwing
away a subset of codimension 2, we may assume that $T$ is locally
free (of course now $X$ need not be compact any longer) and also
that $Z = \emptyset.$  \\
 We claim that - analogously to [PS00,2.1] -
$\sO_Y \oplus T^*$ is a subring of $f_*(\sO_X) = \sO_Y \oplus
\sE^*$ outside a set of codimension $2$. In fact, consider the
multiplication map $$ \sigma: (\sO_Y \oplus T^*) \otimes (\sO_Y
\oplus T^*) \to \sO_Y \oplus \sE^*.$$ In combination with the
projection $\sE^* \to \sF^*,$ we obtain a map $$ \tau: (\sO_Y
\oplus T^*) \otimes (\sO_Y \oplus T^*) \to \sF^*.$$ Since $\sF_C $
is the unique maximal ample subbundle of $\sE_C$ (see [PS00,2.3])
and since $T_C$ is numerically flat, we conclude that $\tau \vert
C = 0.$ Thus $\tau = 0$ generically, hence $\tau = 0.$ Therefore
$\sigma $ maps into $\sO_Y \oplus T^*$, so that $\sO_Y \oplus T^*$
is a subring of $f_*(\sO_X).$ \\
 As in the proof of [PS00,2.1], we put
$Z = {\rm Specan} (\sO_Y \oplus T^*),$ a normal complex space.
Removing a set of codimension $2$, we may assume that $Z$ is
smooth, hence $Z \to Y$ is flat. Since $Z_C \to C$ is \'etale, it
follows  that $Z \to Y$ is unramified, hence \'etale.  \\ By
construction, $\sV^*$ is a subsheaf of $T^*,$ and since $\sV^*_C$
is flat, $\sV^*$ must be a subbundle of $T^*$ outside codimension
$2$. \\
Coming to the original notations, we have now constructed a finite \'etale map $h: Z \to Y \setminus W$
with $W$ in $Y$ being of codimension at least 2.  Since $\pi_1(Y \setminus W) = \pi_1(Y),$ clearly $h$
extends to a covering $h: \tilde Z \to Y$ with $\tilde B$ compactifying $B.$ Furthermore
the map $g: X \setminus f^{-1}(W) \to B$ extends to a map $g: X \to \tilde Z$ such that $f = h \circ g$
everywhere. Finally by construction, $\sV$ is a quotient of $\sE_h.$

\qed

\begin{remark} {\rm Actually the conclusions of
(3.8) hold in a more general setting. Instead of assuming the
existence of $\sV,$ it suffices to assume that $\sE_C$ is not
ample for general $C $ cut out by general elements in $\vert m_i
H_i \vert,$ where $m_i \gg 0,$ so that the theorem of
Mehta-Ramanathan applies. }
\end{remark}

The following generalization of 3.6 now holds:

\begin{corollary} Let $d = 3,$ $\dim Y \geq 3$ and
$h^{1,1}(Y) = 1.$ Then $q(X) = q(Y)$ unless $f$ is \'etale.
\end{corollary}

\proof We start as in the proof of 3.6, but we can only conclude
$\sL' \equiv 0.$ Now we apply 3.8(2) to conclude. \qed

Let us come back to the general context. We introduce the
following notation.

\begin{definition} Let $Y_n$ be a projective manifold.
Then we set (for $k \leq n$ and $l \leq n-k$): $${\rm
Lef}^{k,l}(Y) = \{ \alpha \in H^{k,k}(Y) \vert \wedge \alpha:
H^{n-k+l}(Y) \to H^{n+k+l}(Y) {\rm \ is \ onto } \}.$$
\end{definition}

By the Lefschetz theorem $[\omega^k] \in {\rm
Lef}^{k,l}(Y)$ for all K\"ahler classes $[\omega]$ and for all $k,l.$  Moreover it is clear that ${\rm Lef}^{k,l}(Y)$ is non-empty and
Zariski open in $H^{k,k}(Y).$ \\
Let $\sE_r$ be a vector bundle on the projective manifold $X_n.$ Then $c_i(\sE) > 0$ is to say that $c_i(\sE) \cdot Z > 0$ for all
irreducible subvarieties $Z \subset X$ such that $\dim Z = n-i.$ If $c_r(\sE)$ is represented by a positive form, then of course
$c_r(\sE) > 0$, but the converse is not true.

\begin{problem} Describe the structure of ${\rm Lef}^{k,l}(Y)$.
If $\alpha$ is a class represented by a positive $(k,k)-$form, is
it true that $\alpha \in {\rm Lef}^{k,l}(Y)?$
\end{problem}

In the notation of (3.11), 3.2 can be
restated as follows.

\begin{proposition}  Suppose that
$c_{d-1}(\sE_f) \in {\rm Lef}^{d-1,l}(Y).$ Then $H^i(Y) \to H^i(X) $
is onto for $i \leq n+1-d-l$.
\end{proposition}

There is a class of non-ample vector bundles whose top Chern classes are in ${\rm Lef}^{k,l}(Y)$.
These bundles are provided by the theory developed in [dCM00]. First we set up some definitions.

\begin{definition} Let $X_n$ be a projective manifold and $\sL$ a line bundle on $X.$
\begin{enumerate}
\item $\sL$ satisfies the condition $(HL)_k$ (HL for hard Lefschetz) if and only if
$$ \wedge c_1(\sL)^k: H^{n-k}(X) \to H^{n+k}(X) $$
is injective (isomorphic, equivalently).
\item $\sL$ is HL-ample  if $\sL$ fulfills $(HL)_k$ for all $k.$
\end{enumerate}
\end{definition}

HL-ample bundles are called ``lef'' in [dCM00]. That paper establishes the following result:

\begin{theorem} (de Cataldo, Migliorini) Suppose that $\sL$ is semi-ample with associated morphism $\phi.$
Then $\sL$ is HL-ample if and only if $\phi$ is semi-small.
\end{theorem}

Recall that $\phi: W \to Z $ is semi-small
if $$ \dim Z^k + 2k \leq \dim W $$
for all $k \geq 0,$ where $Z^k$ is the set of all $z \in Z$ such that $\dim \phi^{-1}(z)$ has dimension $k.$


Given a vector bundle $\sE$, we let
$$ \zeta = c_1( \sO_{\bP(\sE)}(1)).$$

\begin{definition} A vector bundle $\sE$ is HL-ample if and only if $\sO_{\bP(\sE)}(1) $ is HL-ample, i.e if and only if
$$ \wedge \zeta^k : H^{\dim \bP -k}(\bP(\sE)) \to H^{\dim \bP +k}(\bP(\sE)) $$
is injective for all $k.$
\end{definition}

Now (3.15) yields:

\begin{corollary} Let $\sE$ be a vector bundle such that  $\sO_{\bP(\sE)}(1) $ is semi-ample. Then
$\sE$ is HL-ample if and only if the canonical morphism $\bP(\sE) \to Z $ associated to  $\sO_{\bP(\sE)}(1) $ is
semi-small.
\end{corollary}

The connection to Problem 3.12 is provided by

\begin{proposition}  Let $X_n$ be a projective manifold and $\sE$ a rank $r-$bundle over $X.$ Let $l \geq 0.$
$$ \wedge c_r(\sE): H^{n-r-l}(X) \to H^{n+r-l}(X) $$
is injective if
$$ \wedge \zeta: H^{n+r-2-l}(\bP (\sE)) \to H^{n+r-l}(\bP(\sE)) $$
is injective.
\end{proposition}

\proof [BG71,1.1] \qed

\begin{corollary} Let $\sE$ be a vector bundle of rank $r$ such that
$\sO_{\bP(\sE)}(1)$ is semi-ample and that
the associated map $\bP(\sE) \to Z$ is semi-small, i.e., $\sE$ is $HL-$ample. Then
$\wedge c_r(\sE): H^{n-r-l}(X) \to H^{n+r-l}(X) $ is an isomorphism for all $l \geq 0.$
\end{corollary}

\proof We verify the injectivity of $\zeta:  H^{n+r-2-l}(\bP (\sE)) \to H^{n+r-l}(\bP(\sE))$ by
realizing $\zeta $ as part of the map
$$ (\wedge \zeta)^{l+2}: H^{n+r-2-l}(\bP (\sE)) \to H^{n+r-l}(\bP(\sE))$$
which is an isomorphism by our assumption. \qed

\begin{corollary} Let $f: X \to Y$ be a $d-$sheeted covering of $n-$dimensional projective manifolds. Suppose that
$\sE_f$ is HL-ample. Then
$$ H^{i}(Y) \to H^{i}(X) $$
is surjective for $i \leq n+1-d.$
\end{corollary}

\proof We verify the condition given in (3.2)
by proceeding as in [So78,1.17] by virtue of (3.18); notice that it is of course irrelevant whether we shift by $-l$ or by $l.$
\qed

Observe that the bundle $\sE$ from example 2.1 is HL-ample! Thus (3.2) gives a new proof of (2.7).

\begin{remark}{\rm  It is also very interesting to consider only
$$ \wedge c_r(\sE): H^{n-r}(X) \to H^{n+r}(X)$$ and also to see what the positivity of $c_r(\sE)$ gives.
In order to see wether this map is injective we do not need the full strength of $(HL)$ for $\zeta;$
we just need $(HL)_1$ - besides spannedness - apart from the fact that the injectivity of $\wedge c_r(\sE)$ is
weaker than the injectivity of $\wedge \zeta.$ Also note that the injectivity of $\wedge c_r(\sE)$ is indeed weaker than the
semi-smallness of $\phi_{\zeta};$ consider e.g., $\sE = T_{\bP_n}(-1).$
A necessary condition is given by the next proposition.}
\end{remark}

\begin{proposition} Suppose that $c_r(\sE_r) > 0$  and that $\sO_{\bP(\sE)}(1)$ is semi-ample.
Let $\phi: \bP(\sE) \to Z $ be the associated map. Then $2\dim T - \dim \bP(\sE) \leq \dim \phi(T)$ for all $T$ with $\dim T
\leq n.$
\end{proposition}

\proof First observe that $c_r(\sE) > 0$ implies $\dim \phi^{-1}(z) \leq r-1$ for all $z \in Z.$ In fact, otherwise there exists an
irreducible $S \subset \bP(\sE) $ of dimension $r$ such that $\dim \phi (S) = 0.$ Now let $\pi: \bP(\sE) \to X$ be the projection; then
$ \pi \vert S$ is finite, so that $S$ is a multi-section over $S' = \pi(S).$ After a base change $h:S \to S',$ the variety $S$ will be
a section. Now $h^*(\sE \vert S') $ must a trivial quotient, hence $c_r(\sE \vert S') = 0.$ This is only possible if
$\dim S' < r,$ contradiction. \\
So $\dim \phi^{-1}(z) - \dim \phi(T) \leq r-1$ and
thus $\dim T - \dim \phi(T) \leq r-1.$ Now use $\dim T \leq n$ to verify that $2\dim T - \dim \bP(\sE) \leq \dim \phi(T)$.
\qed

\section{The structure of the ramification divisor}

Again let us fix a covering $f: X \to Y$ of projective manifolds. Let $n = \dim X = \dim Y$ and $d$ be the degree of $f.$
Let $R \subset X$ be the ramification divisor of $f$ so that $\omega_X = f^*(\omega_Y) \otimes \sO_X(R)$. We assume always $R \ne \emptyset.$
Finally let $B \subset Y$ be the
branch locus.
A standard technique to
obtain information on $R$ in case $Y$ is Fano of index $m$ is to use adjunction theory, applied to (now we use the language of divisors)
$K_X + mL$, where $L = f^*(H)$ where $-K_Y = mH$ with $m$ the index of $Y.$  In particular:

\begin{theorem}  Let $Y$ be projective space or a quadric (of dimension at least 3). Then $R$ is (very) ample.

\end{theorem}

\proof [Ei82] and [LPS89]. \qed

If $Y$ is del Pezzo, then $R = K_X+(n-1)L$. Now adjunction theory says (see [BS95]) that $K_X+(n-1)L$ is big and nef unless $X$ is in a
small list (X is a quadric fibration over a curve or a scroll over a surface), the maps always given by a multiple of $K_X+(n-1)L.$ These special
cases need to be considered separately. Apart from that $K_X+(n-1)L$ is big and nef, hence a multiple is spanned giving a birational
map which is well understood: the result is the so-called first reduction $X'$ which needs to be studied further; actually the first reduction
only contracts $\bP_{n-1}$ with normal bundle $\sO(-1).$
In all other cases
$R$ is automatically ample. It is interesting to notice that in Example 2.1, $R$ is not ample. In fact, we have a $\bP_1-$bundle $q: X \to \bP_2$
and $R = q^*(\sO(2)).$ We will come back to this situation in (4.4)\\
\\
Here we propose another method, namely we use the surjectivity of the canonical map
$$ \kappa: f^*f_*(\omega_{X/Y}) \to \omega_{X/Y} = \sO_X(R)$$
Since $ f^*f_*(\omega_{X/Y}) = \sO_X \oplus f^*(\sE),$ we can use properties of $\sE $ to get information on $R.$

\begin{theorem}
\begin{enumerate}
\item If $\sE$ is spanned resp. nef, then $R$ is spanned resp. nef.
\item If $\sE$ is ample, then $R$ is big and nef. If moreover $R = f^{-1}(B)$ set-theoretically,
or if $-K_Y$ is nef with $b_2(Y) = 1$ or
if some symmetric power $S^p\sE$ is
spanned with $b_2(Y) = 1$, then $R$ is ample.
\end{enumerate}
\end{theorem}

\proof (1) is clear by what we said before. In order to prove (2) we observe that the map $\sO_X \to \sO_X(R)$ induced by $\kappa$ is the
canonical inclusion. Therefore we obtain a surjective map
$$ f^*(\sE) \vert R \to \sO_R(R)$$
so that $\sO_R(R)$ is ample. Hence $R$ itself is big and nef and the only potential curves $C$ preventing $R$ from being ample
are curves $C$ disjoint from $R$ with $R \cdot C = 0.$ \\
If now $R = f^{-1}(B)$, then $B \cdot f(C) = 0.$ However $B$ is ample by [PS00,1.15], since $\det \sE$ is ample.\\
So suppose now that $-K_Y$ is nef. We are going to apply the base point free theorem in order to prove that some multiple
$mR$ is spanned. To do so, we will check that $rR - K_X$ is nef for some positive integer $r$. And indeed, $rR - K_X =
(r-1)R - f^*K_Y$. Now let $\phi: X \to W$ be the associated morphism. Then the exceptional locus $E$ lies completely in the
``\'etale part'' of $f$ and therefore $\phi$ induces a morphism $\psi: Y \to Z$ contracting in particular $f(C).$ This
contradicts $b_2(Y) = 1.$ \\
The last statement follows in the same way, since now $R$ is semi-ample from the beginning. \qed

From (4.2) we can again derive (4.1), but much more can be said:

\begin{corollary}
\begin{enumerate}
\item If $Y$ is a Grassmannian, then $R$ is ample and spanned.
\item If $Y$ is abelian or a product of an abelian variety and a projective space (say), then $K_X = R$ is nef.
\item If $Y$ has $K_Y \equiv 0$ and if $\sE_f$ is nef, then $K_X = R$ is nef. If $\sE$ is ample, then $K_X$ is big
and nef and even ample, if $b_2(Y) = 1.$
\end{enumerate}
\end{corollary}

In fact, $\sE$ is ample in case (1) by [Ma97] and nef in case (2) by [PS00]. Notice that if $A$ is abelian, then $K_X = R$
so that the corollary says that $K_X$ is nef. Of course this follows also from Mori theory: if $K_X$ were not nef, then
there is a rational curve $C \subset X$ (with $K_X \cdot C < 0),$ so that $X$ cannot be a covering of an abelian variety.

\begin{corollary} Let $Y$ be the Fano 3-fold $V_5$ (example 2.1). Then $R$ is ample unless $X$ arises by base change from Example
2.1. To be more precise, let $f_0: X_0 \to Y$ be the triple cover from (2.1) with projection $p: X_0 \to T = \bP_2$.
Then there exists a finite cover $g: S \to T$ such that $X = X_0 \times_T S$ and $f = f_0 \circ q$ with $q: X \to X_0$ the
projection.
\end{corollary}

\proof Since $\sE$ is spanned, $R$ is spanned. Moreover we deduce from 2.4 and 4.2 (proof) that $R.C > 0$ for all curves with the following
(possible) exceptions: $f(C) $ is a line or $C \cap R = \emptyset.$ Let $L = f^*(H),$ with $H$ the ample generator on $Y.$ Then
we have $K_X + 2L = R.$ Suppose that $R$ is not ample. Then $K_X + 2L$ is nef but not ample, hence by adjunction theory (and $R \ne \emptyset$)
we are in one of the following case [BS95]:
\begin{enumerate}
\item $K_X + 2L$ is big and defines a birational map $\phi: X \to X'$
\item $K_X + 2L$ defines a quadric fibration over a curve
\item $K_X + 2L$ defines a scroll $X \to S.$
\end{enumerate}
Let us rule out (1) and (2). In case (1) $\phi$ is just the blow-up of smooth points, so that $\phi$ contracts $E_i = \bP_2$ with
normal bundle $\sO(-1).$ If none of the $E_i$ meets $R$, then we obtain a contradiction to $b_2(Y) = 1$ as in the proof of (4.2).
So suppose that $E_1 \cap R \ne \emptyset.$ Then we conclude that every curve in $E_1$ must be mapped to a line in $Y$ which is
absurd. So (1) is ruled out. (2) is done in exactly the same way (consider a general fiber). \\
So we end up with a scroll $r: X \to S.$ Then $L \cdot r^{-1}(s) = 1$ and $f$ has degree 1 on every $r-$fiber. In particular $f$ maps
the $r-$fibers to lines in $Y.$ Since $T$ is the parameter space of lines in $Y$ and $X$ is the total space of the family, we obtain
a finite map $S \to T$ and $X$ arises from $X_0$ by base change. \qed

It is also interesting to ask, what consequence the ampleness of $R$ for the bundle $\sE$ might have. At least one can state

\begin{proposition} If $R$ is ample, then $\det \sE$ is ample.
\end{proposition}

\proof By [PS00,1.15] we have $2\det \sE = f_*(R) $ (as divisors). Introduce a norm N on $H^{2n-2}(X).$ Via the injective map $f^*(Y):
H^{2n-2}(Y) \to H^{2n-2}(X),$ we obtain a norm on $H^{2n-2}(Y).$ Then there is a positive constant $\epsilon$ such that
$ \vert R \cdot B \vert \geq \epsilon N(B)$ for every effective curve $B \subset X$ since $R$ is positive on the closure of effective
curves.
Via the projection formula we obtain for every irreducible curve $C \subset Y$ that
$$ \vert  2c_1(\sE) \cdot C \vert  = \vert  R \cdot f^*(C) \vert \geq \epsilon N(C)$$
Kleiman's criterion now gives the ampleness of $\det \sE.$ \qed

\section{Further problems and remarks}

In the following we collect some problems on coverings $f: X \to Y$ of projective manifolds of dimension $n.$
Mostly $Y$ will be Fano with $b_2 = 1.$ We do not touch any more problems concerning the topology of low degree coverings, this being
intensively studied in section 3.

\begin{problem} Suppose that $Y$ is Fano with $b_2(Y) = 1$ or even just a projective manifold with $b_2(Y) = 1.$
Is then $\sE_f$ nef? If $H$ denotes the ample generator, is then
$\sE$ spanned outside the base locus of $\vert H \vert?$ What can be said about positivity of Chern classes of $\sE_f$?
\end{problem}

Note that one cannot expect in general $\sE$ to be spanned: let $Y$ be a del Pezzo 3-fold of degree 1, so that $H^3 = 1.$ Then it is known that
$\vert H \vert$ has exactly one base point $y_0$ and $\vert 2H \vert $ is spanned. Now take a smooth member $B \in \vert 2H \vert$ and let
$f: X \to Y$ be the cyclic cover of degree 2 branched along $B.$ Then $\sE_f = H$ and therefore $\sE$ is not spanned at $y_0.$

\begin{problem}Let $Y$ be a Fano manifold with $b_2(Y) = 1.$
Suppose that $\sE$ is ample on all rational curves with minimal degree or ample on all covering families of rational curves of minimal degree
among covering families or on all rational curves $C$ with $-K_X \cdot C \leq \dim Y + 1.$
Is then $\sE$ ample?
\end{problem}

The only evidence for a positive answer is provided by (2.1) and the general philosophy that the geometry of $Y$ is very much dictated
by the minimal (covering) rational curves. A slightly weaker formulation would be:

\begin{problem} Is $\sE_f$ generically ample in the sense that $\sE_f$ is ample say on sufficiently general complete intersections given by
very ample divisors?
\end{problem}

\begin{problem} Fix the Fano manifold $Y$ and assume that $\sE_f$ is ample for all coverings $f: X \to Y.$ Is then $Y$ rational homogeneous?
\end{problem}

Note that necessarily $b_2(Y) = 1,$ since otherwise by Mori theory $Y$ admits a non=trivial map to a normal projective variety.
Again the only support for a positive answer is provided by (2.1). Even for threefolds the answer is not clear. However here an explicit list exists
and one can treat the problem case by case, e.g., the 3-dimensional cubic in $\bP_4.$ \\
It might be interesting to investigate also coverings over other manifolds such as K3 surfaces or Calabi-Yau manifolds.

\begin{problem} Let $Y$ be a K3 surface or a Calabi-Yau manifold. Is $\sE$ always nef?
\end{problem}

Note that if $f:X \to Y$ is a cyclic cover over the K3 surface $Y$, then $\sE$ is spanned. In fact, let $B \subset Y$ be the smooth branch locus;
we may assume $D$ connected, hence $D$ is a smooth irreducible curve. Since $D$ has a root as line bundle, $D$ cannot be a $(-2)-$curve. Hence
$g(D) \geq 1$ and if $g(D) = 1,$ then $D$ is a fiber of an elliptic fibration (apply Riemann-Roch to see $h^0(\sO_Y(D)) \geq 2)$.
Finally if $g(D) \geq 2,$ equivalently $D^2 > 0,$ then the normal bundle sequence reads
$$ 0 \to \sO_Y \to \sO_Y(D) \to N_D = K_D \to 0$$
and therefore $H^1(\sO_Y) = 0$ yields the spannedness of $D.$
Since $\sE_f = \sO_Y(D) \oplus \ldots \oplus \sO((d-1)D),$
this bundle is spanned, too.

\begin{remark} {\rm
 There is an analogy with positive-dimensional fibrations which we
would like to mention. One standard reference here is [Mo87]. So
let $f: X \to Y $ be an algebraic fiber space, say with $X$ and
$Y$ smooth (projective). One says $f$ has snc-branching if the set
(in $Y$) of singular fibers is contained is contained in a divisor
in $Y$ with simple normal crossings. If $f$ has snc-branching,
then all sheaves $R^jf_*(\omega_{X/Y})$ are locally free, and if
moreover $f$ is semi-stable in codimension 1, which is to say that
in codimension 1 all fiber are reduced with only normal crossings,
then these sheaves are even nef (Viehweg). If in addition $\sL$ is
a very ample line bundle on $X$, then Koll\'ar has shown that
$\sL^k \otimes \omega_Y \otimes f_*(\omega_{X/Y}) $ is spanned for
$k > m := \dim Y.$ In particular $f_*(\omega_{X/\bP_m})$
is spanned. \\
Using Mumford regularity, one can e.g., prove  for $Y = \bP_m $ and all $j:$ \\
(1) $R^jf_*(\omega_{X/Y}) $ is spanned outside the singular set of this sheaf. \\
Using semi-stable reduction it also can be seen that \\
(2) $R^jf_*(\omega_{X/Y})(-1)$ is spanned outside the singular set,
unless $H^0(R^{k-j}(\sO_X)^{**}) \ne 0,$ where $k = \dim X - \dim Y.$ \\
In general, $R^jf_*(\omega_{X/Y})$ is nef on every curve which is
not in the closure of the set of all $y \in Y$ such that the fiber
$X_y$ is singular, and not reduced or not with normal
crossings.This is the analogue to Lazarsfeld's theorem that
$\sE_f$ is nef outside the branch locus for coverings $f$.}

\end{remark}

\vspace{1cm}
\small
\begin{tabular}{lcl}
  Thomas Peternell&                  & Andrew J. Sommese \\
  Mathematisches Institut &          & Department of Mathematics\\
  Universit\" at Bayreuth &          & University of Notre Dame\\
  D-95440 Bayreuth, Germany&         &Notre Dame, Indiana 46556-4618, U.S,A,\\
  fax: Germany + 921--552785&        & fax: U.S.A. + 574--631-6579 \\
  thomas.peternell@uni-bayreuth.de & & sommese@nd.edu\\
                                    &&URL: {\tt www.nd.edu/$\sim$sommese}\\
\end{tabular}


\begin{thebibliography}{999999}

\bibitem[An85]{An85} Ando,T.: On extremal rays of higher-dimensional varieties. Inv. Math. 81, 347-357 (1985)

\bibitem[BS95]{BS95} Beltrametti,M.; Sommese,A.J.: The adjunction
theory of complex projective varieties. de Gruyter 1995

\bibitem[BG71]{BG71} Bloch,S.; Gieseker,D.: The positivity of the Chern classes of an
ample vector bundle. Inv. Math. 12, 112-117 (1971)


\bibitem[dCM00]{dCM00} de Cataldo,M.; Migliorini,L.: The hard Lefschetz theorem and
the topology of semismall maps. math.AG/0006187

\bibitem[Ei82]{Ei82} Ein,L.: The ramification divisor for branched coverings of $\bP_k^n.$ Math. Ann. 261,
483-485 (1982)

\bibitem[FN89]{FN89} Furushima,M.; Nakayama,N.: The family of
lines on the Fano threefold $V_5.$ Nagoya Math. J. 116, 111-122
(1989)

\bibitem[Fu90]{Fu90} Fujita,T.: Classification theories of
polarized varieties. London Math. Soc. Lect. Notes Ser. 155 (1990)

\bibitem[Is77]{Is77} Iskovskih,V.A.: Fano 3-folds I. Math.
USSR Izv. 11, 485-527 (1977)

\bibitem[Ki96]{Ki96} Kim,M.: Barth-Lefschetz theorem for branched coverings of Grassmannians. J. reine u. angew.
Math. 470, 109-122 (1996)

\bibitem[La80]{La80} Lazarsfeld,R.: A Barth type theorem for
branched coverings of projective spaces. Math. Ann. 249, 153--162
(1980)

\bibitem[La00]{La00} Lazarsfeld,R.: Nefness modulo the branch locus of the bundle associated to a branched
covering. Appendix to [PS00]. Comm. in Algebra 28, 5598-5599 (2000)

\bibitem[LPS89]{LPS89} Lanteri,A.: Palleschi,M.: Sommese.A.J.: Very ampleness of $K_X \otimes \sL^{\dim X}$
for ample and spanned line bundles $\sL$. Osaka J. Math. 26, 647-664 (1989)

\bibitem[Ma97]{Ma97} Manivel,L.: Vanishing theorems for ample vector bundles. Inv. Math. 127, 401-416 (1997)

\bibitem[Mo87]{Mo87} Mori,S.: Classification of higher-dimensional varieties. Proc. Symp. Pure Math. 46,
269-331 (1987)

\bibitem[PS00]{PS00} Peternell,T.; Sommese,A.J.: Ample vector
 bundles and branched coverings. Comm. in Algebra 28, 5573--5597 (2000)

\bibitem[Re95]{Re95} Reid,M.: Non-normal del Pezzo surfaces. Publ. RIMS 30,
 695-727 (1995)

\bibitem[So78]{So78} Sommese,A.: Submanifolds of abelian varieties. Math. Ann. 233,
229-256 (1978)


\end{thebibliography}
\end{document}